\documentclass[11pt]{article}
\usepackage{amssymb,latexsym,start2e,pstcol,pst-plot,amsthm,amsmath}

\newtheorem{thm}{Theorem}[section]

\newtheorem{question}[thm]{Question}

\begin{document}
\pagestyle{plain}

\title{Pattern Avoidance in Set Partitions
}
\author{
Bruce E. Sagan\\
Department of Mathematics\\
Michigan State University\\
East Lansing, MI 48824-1027\\
USA\\
\texttt{sagan@math.msu.edu}
}

\date{\today\\[10pt]
	\begin{flushleft}
	                                       \\[5pt]
	\end{flushleft}}

\maketitle

\begin{abstract}

The study of patterns in permutations in a very active area of current
research.  Klazar defined and studied an analogous notion of  pattern for set
partitions.  We continue this work, finding exact formulas for the
number of set partitions which avoid certain specific patterns.  In
particular, we enumerate and characterize those partitions avoiding
any partition of a 3-element set.  This allows us to conclude that the
corresponding sequences are P-recursive.  Finally, we define a second
notion of pattern in a set partition, based on its restricted growth
function.  Related results are obtained for this new definition.
\end{abstract}

\pagebreak

\section{Introduction}

The study of patterns in permutations has been very active of late;
see the article of Wilf~\cite{wil:pp} for a survey.
Klazar~\cite{kla:aas,kla:cpsI,kla:cpsII} defined and investigated what
it means for a set partition to 
avoid a pattern, generalizing Kreweras's much-studied notion of a
noncrossing partition~\cite{kre:pnc}.  
Recently, Klazar and Marcus~\cite{km:elb} proved a generalization of
the Marcus-Tardos Theorem~\cite{mt:epm} (which itself provided a
demonstration of the F\"uredi-Hajnal and Stanley-Wilf Conjectures),
that in particular gives the asymptotic growth rate of the number of
set partitions avoiding a given pattern.  Here, our focus will be on exact
enumeration.  To make things precise, we will need some definitions.

Let $\bbP$ and $\bbN$ denote the positive and nonnegative integers,
respectively. 
For $m,n\in\bbN$ we have the interval $[m,n]=\{m,m+1,\ldots,n\}$
with special case $[n]=[1,n]$.
If $S$ is any set, then a {\it partition $\pi$ of $S$\/} is a set of
nonempty subsets $B_1,B_2,\ldots,B_k$ of $S$ such that 
$\uplus_i B_i =S$ (disjoint union).  We will write $\pi\ptn S$ and
$\pi=B_1/B_2/\ldots/B_k$.  The subsets are called {\it blocks\/} and
the number of blocks will be denoted
$b(\pi)$.  Most often we will also not use set braces and commas in
the blocks unless they are needed for clarity.
For example, if $\pi=14/2/356$ then $\pi\ptn[6]$ and $b(\pi)=3$.
We will use the notation
\beq
\label{Pin}
\Pi_n=\{\pi\ :\ \pi\ptn[n]\}\qmq{and} \Pi=\biguplus_{n\ge0} \Pi_n.
\eeq

In order to connect Klazar's definition of pattern with the usual one
for permutations, it is convenient to introduce a standardization
map.  If $S$ is any subset of the integers with cardinality $\#S=n$
then the corresponding {\it standardization map\/} is the unique
order-preserving bijection $\sta_S:S\ra [n]$.  When $S$ is clear from
context, we will drop the subscript.  We let $\sta_S$ act element-wise
on objects built using $S$ as label set.  For example, if
$S=\{3,4,6\}$ then $\sta(3)=1$, $\sta(4)=2$, $\sta(6)=3$.
Consequently, for the sequence $p=4346$ we have $\sta(p)=2123$ and for
the partition $\pi=36/4$ we have $\sta(\pi)=13/2$.  The definition of
pattern containment for permutations can now be stated as follows:
If $p=a_1a_2\ldots a_r$ and $q=b_1b_2\ldots b_s$ are permutations, then
{\it $q$ contains $p$ as a pattern\/} if there is a subsequence
$q'=b_{i_1}b_{i_2}\ldots b_{i_r}$ of $q$ with $\sta(q')=p$.  
Otherwise $q$ {\it avoids\/} $p$.  Given a
pattern permutation $p$, we let
$$
\fS_n(p)=\{q\in\fS_n\ :\ \mbox{$q$ avoids $p$}\},
$$
where $\fS_n$ is the symmetric group on $[n]$.
We will also let $\fS=\uplus_{n\ge0} \fS_n$.

For pattern containment in set partitions, we will need the notion of
a subpartition.  A {\it subpartition\/} of $\si$ is
a partition $\si'$ such that each block of $\si'$ is contained in a
different block of $\si$.  For example, $\si=14/236/5$ has $\si'=26/4$
as a subpartition, but not $26/3$ since both $26$ and $3$ are in the
same block of $\si$.  If $\pi$ and $\si$ are set partitions, then 
{\it $\si$ contains $\pi$ as a pattern\/} if
there is a subpartition $\si'$ of $\si$ with $\sta(\si')=\pi$.  Also,
$\si'$ is called a {\it copy\/} of $\pi$ in $\si$.  If $\si$ has no
copies of $\pi$ then it {\it avoids\/} $\pi$.  Continuing our example,
$\si=14/2/356$ contains four copies of the pattern $13/2$, namely
$14/2,\ 14/3,\ 35/4$, and $36/4$.  On the other hand,
$\si$ avoids $134/2$ since any copy of this pattern would
have to have $356$ as a block and then there is no integer that can
take the place of the $2$ in the pattern.
Parallel to the notation above, given a pattern $\pi$ we let
$$
\Pi_n(\pi)=\{\si\in\Pi_n\ :\ \mbox{$\si$ avoids $\pi$}\}
\qmq{and} \Pi(\pi)=\biguplus_{n\ge0} \Pi_n(\pi).
$$
Note that the noncrossing partitions may be defined as those in
$\Pi(13/24)$.

In the following section we will provide exact formulas and generating
functions for $\#\Pi_n(\pi)$
for various patterns $\pi$, including all $\pi\ptn[3]$.
Gessel~\cite{ges:sfp} and Noonan-Zeilberger~\cite{nz:epp} initiated
the study of P-recursiveness 
and its relationship to patterns in permutations.  In
section~\ref{Pre} we consider analogous results for set partitions.
The section after that uses restricted growth functions to give a second
definition of  pattern in a set partition, and various results using
this new notion are presented.  We end with a section
outlining future work and open problems.

\section{Enumeration}
\label{enu}

As is often the case when dealing with set partitions, exponential
generating functions will be useful.  So we begin by setting up some
notation for them.

If $I$ is a set of nonnegative integers, then let
\beq
\label{FI}
F_I(x)=\sum_{i\in I} \frac{x^i}{i!}.
\eeq
We will also use the following notation for a special case
of~\ree{FI} which will appear repeatedly
$$
\exp_m(x)=\sum_{n=0}^{m} \frac{x^n}{n!}.
$$
The next result follows from standard manipulation of 
exponential generating functions 
(see Wilf's book~\cite[Chapter 3]{wil:gen}), so we omit the proof.
\bpr
\label{egf}
Let
$$
a_{n,l}^I=\#\{\si=C_1/C_2/\ldots/C_l\in\Pi_n\ :\ 
\mbox{$\#C_j\in I$ for $1\le j\le l$}\}.
$$
It follows that
$$
\sum_{n=0}^{\infty} a_{n,l}^I \frac{x^n}{n!} = \frac{F_I(x)^l}{l!}.\qqed
$$ 
\epr
Finally, given a pattern $\pi$, we let
$$
F_\pi(x)=\sum_{n=0}^{\infty} \#\Pi_n(\pi) \frac{x^n}{n!}.
$$
This will cause no confusion with~\ree{FI}, since the context will
always make it clear whether the subscript refers to an index set or a
partition.  

For our first result, we will consider the extreme cases where $\pi$
has only singleton blocks or is itself a single block.  When
$\pi=1/2/\ldots/m$, $\si$ contains a copy of $\pi$ if and only if
$\si$ has at least $m$ blocks from which to take the singletons.
Similarly, if $\pi=12\ldots m$ then a copy of $\pi$ can come from any
block of $\si$ having at least $m$ elements.  Combining these
observations with the previous proposition proves the following.
\bth
\label{01}
We have
\bea
\Pi(1/2/\ldots/m)&=&\{\si\in\Pi\ :\ b(\si)<m\},\\[5pt]
F_{1/2/\ldots/m}(x)&=&\exp_{m-1}(\exp(x)-1),\\[5pt]
\Pi(12\ldots m)&=&\{\si\in\Pi\ :\ \mbox{$\#C<m$ for all $C\in\si$}\},\\[5pt]
F_{12\ldots m}(x)&=&\exp(\exp_{m-1}(x)-1).\qqed
\eea
\eth

In what follows, we will often abbreviate
$$
\zh_m=1/2/\ldots/m\qmq{and}\oh_m=12\ldots m.
$$
This is because these elements are the unique minimum and maximum of
the partition lattice.

We can characterize the set partitions which avoid another infinite
families of patterns.
Suppose that $\si=C_1/C_2/\ldots/C_l\ptn S$ and $T\sbe S$.   Then the 
{\it  restriction of $\si$ to $T$\/} is the partition $\si_T$ whose
blocks are the nonempty sets of 
the form $C_i\cap T$, $1\le i\le l$.  Using $\si=14/2/356$ as usual
and $T=\{3,4,6\}$, we obtain $\si_T=36/4$.  Note that $\si'$ is a
subpartition of $\si$ if and only if $\si'=\si_T$ for some $T$.  If
$\si\ptn[n]$ then we 
will use the abbreviations $\si_{\le k}$ and $\si_{>k}$ for the
cases when $T=[k]$ and $T=[k+1,n]$, respectively.  In the following
theorem, we use the falling factorial notation
$$
\fall{k}{i}=k(k-1)\cdots(k-i+1)
$$
as well as the $a_{n,l}^I$ as defined in Proposition~\ref{egf}.
\bth
\label{01+}
We have
\begin{eqnarray}
\label{12}
\Pi(12/3/4/\ldots/m)&=&\{\si\in\Pi\ :\
\mbox{$\exists k$ with $\si_{\le k}=\zh_k$ and $b(\si_{>k})<m-1$}\},\\[5pt]
\label{12sum}
\#\Pi_n(12/3/4/\ldots/m)&=&1+
\sum_{k=1}^{n-1}\sum_{j=1}^{m-2} a_{n-k,j}^\bbP\  \sum_{i=1}^j
{j-1\choose i-1} \fall{k}{i}.
\end{eqnarray}
\eth
\pf\
Suppose $\si\ptn[n]$.  If
$\si=\zh_n$ then clearly $\si$ is a member of both sets in~\ree{12}.
So assume $\si\neq\zh_n$.  Define $k$ to be the largest integer such
that all the elements of $[k]$ are minima of their blocks in $\si$.
Then $\si_{\le k}=\zh_k$, and
$k+1$ is in a block of $\si$ which also contains an element $s\le k$.

Suppose first that $\si\in\Pi(12/3/4/\ldots/m)$.    To show that $\si$
must then be in the left-hand side of~\ree{12} we assume, towards a
contradiction, that $b(\si_{>k})\ge m-1$.  But then $\si_{>k}$
contains a copy, $\si'$, of $\zh_{m-1}$ and we can
take that copy to contain $k+1$ since the minima of any $m-1$ blocks
will do.  Inserting the element $s$ into the block of $k+1$ in $\si'$
gives a copy of $12/3/4/\ldots/m$  in $\si$, a contradiction.

For the reverse inclusion, take $\si\neq\zh_n$ in the right-hand set
of~\ree{12}.  We again proceed by contradiction, assuming that $\si$
contains a copy, $\si'$, of $12/3/4/\ldots/m$ .  Let $k'$ be the element of $\si'$
playing the r\^ole of the $2$ in $12/3/4/\ldots/m$.  Then $k'>k$ since the
elements of $[k]$ are all in separate blocks of $\si$.  Thus the
elements of $\si'$ corresponding 
to the elements $[2,m]$ in $12/3/4/\ldots/m$ are
all in $\si_{>k}$, and are also all in separate blocks.  This
contradicts $b(\si_{>k})<m-1$ and finishes the proof of~\ree{12}.

To obtain the count~\ree{12sum}, we enumerate the elements in the
right-hand set of~\ree{12}.  The $1$ in the sum accounts for the
partition $\si=\zh_n$.  Let $k$ be as  defined in the second paragraph of the
proof.  Then $1\le k\le n-1$ since we are now considering
$\si\neq\zh_n$.  Let $j=b(\si_{>k})$ so, by definition of the $\si$
being counted, $1\le j\le m-2$.  Every block of $\si$ if of one
of the three forms $\{s\}$, $C$, or $C\cup\{s\}$ 
where $s\le k$ and $C\in\si_{>k}$.
Let $i$ be the number of blocks of the third type.  So $i\ge1$ since,
by maximality of $k$, the block of $\si$ containing $k+1$ must be of this form.
Also $i\le j$ by their definitions.  Thus we have verified the limits on
the summations in~\ree{12}.

To count the number of $\si$ for given $i,j,k$ we first note that the
choice of $\si_{\le k}$ is unique and there are $a_{n-k,j}^{\bbP}$
choices for $\si_{> k}$.  To determine $\si$ from these two
subpartitions, it suffices to specify the blocks of type three.  We
already know that the block containing $k+1$ must be of this type, so
there are ${j-1\choose i-1}$ ways to choose the rest of the blocks of
$\si_{> k}$ that will be used.  Now these blocks
(including the one containing $k+1$) can each be unioned with a
unique element $s\le k$ in a total of $\fall{k}{i}$ ways.  This gives
the summand in~\ree{12} and completes the proof of this equation and
of the theorem.
\Qqed

Representing a permutation $p$ by its permutation matrix, the dihedral
group of the square acts on $\fS_n$.  The number of permutations
avoiding a pattern is the same for any two patterns in the same
orbit.  More generally, two permutations $p,q$ are called 
{\it Wilf  equivalent\/}  if 
$\#\fS_n(p)=\#\fS_n(q)$ for all $n\ge0$. For example, it is well known
that any two permutations in $\fS_3$ are Wilf equivalent.

Only one of the symmetries for permutations remains for set partitions.
Given $\pi=B_1/\ldots/B_k\ptn[m]$, define its {\it complement\/}
to be the partition $\pi^c=B_1^c/\ldots/B_k^c$ where 
$$
B_i^c=\{m+1-b\ :\ b\in B_i\}
$$
for $1\le i\le k$.  For example, $(14/2/356)^c=63/5/421$.  The proof
of the next result is trivial and so is omitted.
\ble
\label{c}
For any pattern $\pi$, we have
\bea
\Pi_n(\pi^c)&=&\{\si^c\ :\ \si\in\Pi_n(\pi)\},\\
\#\Pi_n(\pi^c)&=&\#\Pi_n(\pi).\qqed
\eea
\ele
We will call partitions $\pi,\si$ {\it Wilf equivalent\/} if 
$\#\Pi_n(\pi)=\#\Pi_n(\si)$ for all $n\ge0$. 
So, for example, the preceding lemma gives us the Wilf equivalence 
$$
\#\Pi_n(1/2/\ldots/m-2/m-1\ m)=\#\Pi_n(12/3/4/\ldots/m)
$$
with the cardinality of the right-hand side being given in~\ree{12sum}.

We will now give a complete characterization and enumeration of
$\Pi_n(\pi)$ for all $\pi\ptn[3]$.  To do so, it will be useful to
have a few more definitions.  Call $\si$  a {\it matching\/} if
$\#C\le2$ for all $C\in\si$. 
Also, define the {\it double factorial\/}
$$
(2i)!!=1\cdot 3\cdot 5\cdots (2i-1),
$$
which is the number of matchings on $2i$ elements where every block
has size two.
Finally, we say that $\si$ it {\it layered\/} if
it has the form
$$
\si=[i,j]/[j+1,k]/[k+1,l]/\ldots/[m+1,n]
$$
for certain $i,j,k,l,\ldots,m,n$.  
\bth
\label{Pi3}
We have
\begin{eqnarray}
\Pi_n(1/2/3)&=&\{\si\in\Pi_n\ :\ b(\si)\le2\},\\[5pt]
\#\Pi_n(1/2/3)&=&2^{n-1},\\[5pt]
\label{123}
\Pi_n(123)&=&\{\si\in\Pi_n\ :\ \mbox{$\si$ is a matching}\},\\[5pt]
\#\Pi_n(123)&=&\sum_{i\ge0} {n\choose 2i}(2i)!!,\\[5pt]
\Pi_n(12/3)&=&\{\si\in\Pi_n\ :\ 
\mbox{$\exists k$ with $\si_{\le k}=\zh_k$ and $\si_{>k}=[k+1,n]$}\},\\[5pt]
\#\Pi_n(12/3)&=&1+ {n\choose 2},\\[5pt]
\Pi_n(1/23)&=&\{\si\in\Pi_n\ :\ 
\mbox{$\exists k$ with $\si_{\le k}=[1,k]$ and $\sta(\si_{>k})=\zh_{n-k}$}\},\\[5pt]
\#\Pi_n(1/23)&=&1+ {n\choose 2},\\[5pt]
\label{13/2}
\Pi_n(13/2)&=&\{\si\in\Pi_n\ :\ \mbox{$\si$ is layered}\},\\[5pt]
\label{n13/2}
\#\Pi_n(13/2)&=&2^{n-1}.
\end{eqnarray}
\eth
\pf\
All these equations except the last two are easy consequences of
Theorems~\ref{01} and~\ref{01+} and Lemma~\ref{c}.

To prove~\ree{13/2}, first note that it is clear from the definition of
``layered'' that such a partition can not have a copy of the pattern $13/2$.
For the reverse direction, suppose $\si$ avoids $13/2$ and let $C$
be the block of $\si$ containing $1$.  Also, let $i=\max C$.  We claim
that $C=[1,i]$.  This is clear if $i=1$.  If $i>1$ then suppose,
towards a contradiction, that there is some $j$ with $1< j< i$ and
$j\not\in C$.  But then $1i/j$ is a copy of $13/2$ in $\si$, a
contradiction.  Considering the block of $\si$ with minimum $i+1$ and
iterating this process completes the proof of~\ree{13/2}.

To prove~\ree{n13/2}, just note that any layered partition can be
obtained from the sequence $12\ldots n$ by inserting slashes in the
$n-1$ spaces between the numbers.  \Qqed

We should note that Klazar also mentioned~\ree{123} 
in~\cite[Example 1]{kla:cpsI}, thus showing that $\pi=123$ is  
a set partition with superexponential growth rate.  (Something which
can not happen for permutations.)

\section{P-recursion}
\label{Pre}

We now use the results of the previous section to investigate
when various sequences of the form $\#\Pi_n(\pi)$, $n\ge0$, are
P-recursive.  A sequence $(a_n)_{n\ge0}$, is {\it P-recursive\/}
(polynomially recursive) if
there are polynomials $P_0(n),P_1(n),\ldots,P_k(n)$ (not all zero)
such that
$$
P_0(n)a_n+P_1(n)a_{n+1}+\cdots+P_k(n)a_{n+k}=0
$$
for all $n\ge0$.
As a simple example, the sequence with elements $a_n=n!$ is P-recursive since we
always have $(n+1)a_n-a_{n+1}=0$.

Gessel~\cite{ges:sfp} first mentioned the problem of determining for
which permutations $p$ the sequence $a_n=\#\fS_n(p)$, $n\ge0$, is  P-recursive.
Noonan and Zeilberger~\cite{nz:epp} conjectured that the sequence is
P-recursive 
for all $p$, although later evidence has caused Zeilberger to change
his mind~\cite{ev:pcp} and conjecture that  it is not P-recursive for $p=1324$.
For set partitions, the numbers $\#\Pi_n(\pi)$ do not always form a P-recursive
sequence, as we will show shortly.  To do so, we need to introduce
some ideas from the theory of D-finite power series.

Let $f(x)$ be a formal power series.  Then $f(x)$ is {\it D-finite\/}
(differentiably finite) if there are polynomials
$p_0(x),p_1(x),\ldots,p_k(x)$ (not all zero)
such that 
\beq
\label{D}
p_0(x)f(x)+p_1(x)f'(x)+\cdots+p_k(x)f^{(k)}(x)=0.
\eeq
A simple example is the function $f(x)=e^x$ which
satisfies $f(x)-f'(x)=0$.  Stanley~\cite{sta:dfp}
was the first to bring the theory of D-finite series, which had long
been used for differential equations, to bear on combinatorial problems.
We will need the following two results of his, the first of which can also be
found in the work of Jungen~\cite{jun:sts}.
\bth[Jungen~\cite{jun:sts}, Stanley~\cite{sta:dfp}]
\label{js}
A sequence $(a_n)_{n\ge0}$ is P-recursive if and only if its ordinary
generating function $f(x)=\sum_{n\ge0} a_n x^n$ is D-finite.\qqed
\eth

\bth[Stanley~\cite{sta:dfp}]
\label{st}
If $(a_n)_{n\ge0}$ and $(b_n)_{n\ge0}$ are P-recursive sequences, then so
is their point-wise product $(a_n b_n)_{n\ge0}$.\qqed
\eth

\bco
\label{Pexp}
A sequence $(a_n)_{n\ge0}$ is P-recursive if and only if its exponential
generating function $F(x)=\sum_{n\ge0} a_n x^n/n!$ is D-finite.
\eco
\pf\
We will only prove the reverse implication, as the forward direction
is obtained from that proof by just reversing the steps.  So suppose
$F(x)$ is D-finite.  Then by Theorem~\ref{js}
the sequence $(a_n/n!)_{n\ge0}$ is P-recursive.  Also, we have already
seen that the sequence $(n!)_{n\ge0}$ is P-recursive.  Thus, by 
Theorem~\ref{st}, $(a_n)_{n\ge0}$ is P-recursive.
\Qqed

For an example where the $\#\Pi_n(\pi)$ do not form a P-recursive sequence,
consider the pattern $\pi=\ep$, the empty partition.  So
$\Pi_n(\ep)=B_n$, the $n$th Bell number. 
\bpr
The sequence $\#\Pi_n(\ep)$, $n\ge0$, is not P-recursive.
\epr
\pf\
Suppose, towards a contradiction, that this sequence is P-recursive.
Using Proposition~\ref{egf}, we get the well-known generating function
for the Bell numbers
$$
F_\ep(x)=e^{e^x-1}.
$$
By the previous corollary, $F_\ep(x)$ must be D-finite and so must
satisfy~\ree{D} for certain polynomials $p_i(x)$.  Taking the
derivatives and dividing by $F_\ep(x)$ which is never zero, we get an
equation of the form
\beq
\label{q}
q_0(x)+q_1(x)e^x+q_2(x)e^{2x}+\cdots+q_k(x)e^{kx}=0
\eeq
where 
$$
q_i(x)=p_i(x)+\sum_{j>i} a_{i,j} p_j(x)
$$
for certain constants $a_{i,j}$.  So since the $p_i(x)$ are
polynomials which are not all zero, the same must be true of the $q_i(x)$.
But this implies that $e^x$ is an algebraic function, a contradiction.
\Qqed

\begin{question}
For what set partitions $\pi$ is the sequence $\Pi_n(\pi)$, $n\ge0$, P-recursive?
\end{question}

We will now show that all of the patterns considered in the previous
section give rise to P-recursive sequences.  To do so, we will need a
few more definitions and results.  In his work on the growth rate of
$\#\Pi_n(\pi)$, Klazar~\cite{kla:cpsI} was lead to consider the
following patterns.  A {\it sufficiently restricted partition\/} or
{\it srp\/} is a matching $\pi$ such that, if $S$ is the union of the
doubletons in $\pi$, then
$$
\sta(\pi_S)=1 a_1/2 a_2/\ldots/k a_k
$$
for some permutation $a_1 a_2\ldots a_k$ of $[k+1,2k]$.
\bth[Klazar~\cite{kla:cpsI}]  
\label{kla}
If $\pi$ is an srp then the ordinary
generating function for the sequence $\#\Pi_n(\pi)$, $n\ge0$, is
rational with integer coefficients.  In particular, this sequence is P-recursive.
\eth

We will also need the following result.
\bth[Stanley~\cite{sta:dfp}]  
\label{sta}
If $f(x)$ is D-finite and $g(x)$ is
algebraic with $g(0)=0$, then the composition $f(g(x))$ is D-finite.
\eth

\bth
\label{Pthm}
For $m\ge1$, the following sequences are P-recursive as $n$ varies
over $\bbN$: 
$$
\#\Pi_n(1/2/\ldots/m),\ \#\Pi_n(12\ldots m), \qmq{and}
\#\Pi_n(12/3/4/\ldots/m).
$$
Furthermore, for any $\pi\ptn[3]$ the sequence $\#\Pi_n(\pi)$, $n\ge0$,
is P-recursive.
\eth
\pf\
The only one of these sequence which is not covered by
Theorem~\ref{kla} is the one for $12\ldots m$.  But in
Theorem~\ref{01} we noted that the exponential generating function for
this pattern is $\exp(\exp_{m-1}(x)-1)$.  We have already seen that
$f(x)=\exp(x)$ is D-finite.  And $g(x)=\exp_{m-1}(x)-1$ is algebraic
since it is a polynomial.  So we are done by Theorem~\ref{sta} and
Corollary~\ref{Pexp}.  
\Qqed

\section{Restricted growth functions}
\label{rgf}

There is a second, natural definition of pattern containment for set
partitions which arises from considering them as restricted growth
functions.  In order to make this connection, we will write all of our
partitions $\pi=B_1/B_2/\ldots/B_k$ in {\it canonical order\/} which
means that the blocks are indexed so that
\beq
\label{can}
\min B_1<\min B_2<\ldots<\min B_k.
\eeq
If $S\sbe B_j$ for some $j$ then it will also be convenient to use the notation
\beq
\label{BS}
B(S)=j.
\eeq

A {\it restricted growth function\/} (RGF) is a sequence 
$r=a_1 a_2\ldots a_n$ of positive integers such that
\ben
\item $a_1$=1, and
\item for $i\ge2$ we have $a_i\le 1+\dil\max_{j< i} a_j$.
\een
The number of elements of $r$ is called the {\it length\/} of $r$ and
denoted $l(r)$.   For example, $r=123133$ is a restricted growth
functions with $l(r)=6$, while $r=123153$ is not an RGF because there
is no $4$ in the prefix before the $5$.  Let
$$
R_n=\{r\ :\ \mbox{$r$ an RGF with $l(r)=n$}\}
\qmq{and} R=\biguplus_{n\ge0} R_n.
$$

There is a well-known bijection  $\rho:\Pi_n\ra R_n$.  Given
$\pi\in\Pi_n$ in canonical order, we let $\rho(\pi)=a_1a_2\ldots a_n$ where
\beq
\label{rho}
a_i=B({i}).
\eeq
If one considers the example partition $\pi=14/2/356=B_1/B_2/B_3$ from the
introduction, then $\rho(\pi)=123133$ (the example RGF above).
The definition~\ree{can} of ``canonical order'' ensures that $\rho(\pi)$
is an RGF.  Furthermore, it is easy to construct an inverse for $\rho$
using~\ree{rho}.  So we can work with a partition or its RGF
interchangeably.  

We now define pattern containment in $R$ analogously to the way it is
defined for permutations.  If $r\in R_k$ is the pattern RGF, then we say that 
$s=b_1 b_2\ldots b_n\in R_n$ {\it contains\/} $r$ if there is a
subsequence $b_{i_1}b_{i_2}\ldots b_{i_k}$ of $s$ which standardizes
to $r$.  Otherwise $s$ {\it avoids\/}  $r$.  By way of illustration, if
$r=121$ then there are two copies of $r$ in 
$b_1 b_2 b_3 b_4 b_5 b_6=123133$, namely $b_1 b_2 b_4=121$ and 
$b_1 b_3 b_4 =131$.

If $\pi$ and $\si$ are such that $\rho(\si)$ contains $\rho(\pi)$ then
we say that $\si$ {\it R-contains\/} $\pi$, and that $\si$
{\it  R-avoids\/} $\pi$ otherwise.  We will also add an ``R'' prefix
to other terms defined 
in the introduction in order to refer to this new definition.
We can see R-containment directly
in terms of partitions as follows:  
$\si$ R-contains $\pi$ if and only if $\si$ has a subpartition
$\si'=C_1'/C_2'/\ldots/C_k'$ (in canonical order) with
$\sta(\si')=\pi$ and
\beq
\label{BC}
B(C_1')<B(C_2')<\ldots<B(C_k').
\eeq
For example, of the four copies of $13/2$ in
$14/2/356$ only two of them, namely $14/2$ and $14/3$, are
R-copies.  In fact, this is just a restatement in terms of partitions of the example at the
end of the previous paragraph.  Given a set partition.
$\pi$ we let
$$
R_n(\pi)=\{s\in\Pi_n\ :\ \mbox{$\si$ R-avoids $\pi$}\}
\qmq{and} R(\pi)=\biguplus_{n\ge0}R_n(\pi).
$$
The next proposition is clear from the definitions.
\bpr
\label{RPi}
For every $\pi\in\Pi$ and every $n\ge0$,
\bea
R_n(\pi)&\supseteq&\Pi_n(\pi),\\
\#R_n(\pi)&\ge&\#\Pi_n(\pi).\qqed
\eea
\epr

Note that if $\pi=\zh_m$ or $\oh_m$ then~\ree{BC}
is automatic.  So the next result follows
immediately from the previous proposition and, because of
Theorem~\ref{01}, the corresponding enumerations have already been
done.
\bth
\label{RPim}
We have
\bea
R(1/2/\ldots/m)&=&\{\si\in\Pi\ :\ b(\si)<m\},\\[5pt]
R(12\ldots m)&=&\{\si\in\Pi\ :\ \mbox{$\#C<m$ for all $C\in\si$}\}.\qqed
\eea
\eth

We now turn to R-avoidance of patterns in $\Pi_3$.  In this context,
complementation does 
not necessarily preserve the number of avoiding partitions.  So it
is somewhat surprising that four of the five elements of $\Pi_3$
are R-Wilf equivalent.
\bth
We have
\begin{eqnarray*}
R_n(1/2/3)&=&\{\si\in\Pi_n\ :\ b(\si)\le2\},\\[5pt]
\#R_n(1/2/3)&=&2^{n-1},\\[5pt]
R_n(123)&=&\{\si\in\Pi_n\ :\ \mbox{$\si$ is a matching}\},\\[5pt]
\#R_n(123)&=&\sum_{i\ge0} {n\choose 2i}(2i)!!,\\[5pt]
R_n(12/3)&=&\{\si\in\Pi_n\ :\ 
\mbox{$\exists k$ with $\si_{\le k}=\zh_k$, $\si_{>k}=D_1/\ldots/D_k$
layered, and $B(D_1)>\ldots>B(D_k)$}\},\\[5pt]
\#R_n(12/3)&=&2^{n-1},\\[5pt]
R_n(1/23)&=&\{\si=C_1/\ldots/C_l\in\Pi_n\ :\ 
\mbox{$\#C_i=1$ for $i\ge2$}\},\\[5pt]
\#R_n(1/23)&=&2^{n-1},\\[5pt]
R_n(13/2)&=&\{\si\in\Pi_n\ :\ \mbox{$\si$ is layered}\},\\[5pt]
\#R_n(13/2)&=&2^{n-1}.
\end{eqnarray*}
\eth
\pf\
The equations involving $1/2/3$ and $123$ follow from
Theorems~\ref{Pi3} and~\ref{RPim}.

Now suppose $\si\in R_n(12/3)$ and let $k$ be the largest integer such that
the elements of $[k]$ are the minima of their blocks in $\si$.  So
either $k=n$ (and $\si=\zh_n$) or $k+1$ is in a block $C_i$ of $\si$
where $i\le k$.  Let $m$ be the maximum of $C_i$.  We claim that
$C_i=\{i\}\uplus[k+1,m]$.  If this were not the case then there would
have to be some $l$ with $k+1<l<m$ and with $l\in C_j$ for $j\neq i$.
If $j<i$ then $jl/m$ is an R-copy of $12/3$ in $\si$, and if $j>i$
then $i,k+1/l$ is such an R-copy.  So in either case we have a
contradiction.  Iterating this argument shows that $\si$ has the form
described in the theorem.  It is also clear that partitions of this
form do not have any R-copies of $12/3$, so this completes the
characterization of such partitions.  

To enumerate $R_n(12/3)$, keep $k$ as in the previous paragraph.  Then
the number of $\si$ for a given $k$ is just the number of ways to
distribute the elements of $\si_{>k}$ among the $k$ blocks.  Since
$\si_{>k}$ is layered, this is equivalent to counting the number of
compositions (ordered integer partitions) of $n-k$ into $k$ parts
where $0$ is allowed as a part.  It is well-known that the number of
such compositions is ${n-1\choose k-1}$.  So the total number of $\si$
is $\sum_k {n-1\choose k-1}=2^{n-1}$.  

Next consider $R_n(1/23)$.  Clearly $\si=C_1/\ldots/C_l$ can not
contain an R-copy of $1/23$ if all
blocks other than $C_1$ are singletons.  And if some block of $\si$
other than $C_1$ contains two elements $i,j$ then $1/ij$ is an R-copy
of $1/23$ in $\si$.  So this gives us the required set equality.  Also,
the number of such $\si$ is just the number of choices for $C_1$, which
is $2^{n-1}$ since we must have $1\in C_1$ for the blocks to be
canonically ordered.

Finally, look at $R_n(13/2)$.  From Theorem~\ref{Pi3} and
Proposition~\ref{RPi} we have that $R_n(13/2)$ contains every layered
permutation.  The proof of the reverse containment is the same as that
given for the corresponding containment in~\ree{13/2}, just noting
that the copy of $13/2$ constructed there 
is, in fact, an R-copy.  Of course, this means that the enumeration is
the same as well.
\Qqed

Since no numerically new sequences have been discussed in this
section, we can use Theorem~\ref{Pthm} to conclude the following.
\bth
For $m\ge1$, the following sequence are P-recursive when $n$ varies
over $\bbN$: 
$$
\#R_n(1/2/\ldots/m) \qmq{and}
\#R_n(12\ldots m).
$$
Furthermore, for any $\pi\ptn[3]$ the sequence $\#R_n(\pi)$, $n\ge0$,
is P-recursive.\qqed
\eth

\section{Open problems and new directions}
\label{fwo}

\subsection{The patterns $12/3/4/\ldots/m$ and $1/23\ldots m$}

We were unable to simplify the summation given for 
$\#\Pi_n(12/3/4/\ldots/m)$.  It would be
interesting to do so, or to use them to find the corresponding
exponential generating function. 

We have given characterizations of $\Pi(\pi)$ where $\pi$ is the
minimum, maximum, or one of the atoms in the partition lattice.  It is also
possible to do so for the coatom $1/23\ldots m$.  We did not mention this
earlier because the description is not used for
any of the other results presented.  But we will give it here in case
it turns out to be useful in later work.  To describe the
$\si=C_1/\ldots/C_l$ in $\Pi(1/23\ldots m)$ we assume, as usual, that
$\si$ is written in canonical form.  We will also need the parameter
$c=c(\si)$ which will be the $(m-1)$st largest element of $C_1$, or
$0$ if $\# C_1<m-1$.  Then
$$
\Pi(1/23\ldots m)=\{\si=C_1/\ldots/C_l\in\Pi\ :\ 
\mbox{$\#C_i<m-1$ for $i\ge2$ and $\min C_2>c(\si)$}\}.
$$
The proof of this equality is much like the one for~\ree{12}, where
the first restriction on $\si$ ensures that there can be no copy of
the pattern where the subset corresponding to $23\ldots m$ is in a
block of index at least two, and the second restriction does the same
for $C_1$.

\subsection{Wilf equivalence}

As previously mentioned, any two permutations in $\fS_3$ are Wilf
equivalent.  Babson and West~\cite{bw:p13} showed that 
$123a_3\ldots a_n$ and $321a_3\ldots a_n$ are Wilf equivalent for any
permutation $a_3\ldots a_n$ of $[4,n]$.  This work was later
generalized by Backelin, West, and Xin~\cite{bwx:wes}.  Are the Wilf
equivalences that appeared for both containment and R-containment
isolated incidents or part of a larger picture?

\subsection{Multiple restrictions}

Let $\cP\sbe\fS$ be any set of permutations and define
$$
\fS_n(\cP)=\{q\in\fS_n\ :\ \mbox{$q$ avoids $p$ for all $p\in\cP$}\}.
$$
Simion and Schmidt~\cite{ss:rp} enumerated all such sets where
$\cP\sbe\fS_3$.  Similarly, for $\cP\sbe\Pi$ one can let
$$
\Pi_n(\cP)\ =\ \{\si\in\Pi_n\ :\ \mbox{$\si$ avoids $\pi$ for all $\pi\in\cP$}\}.
$$
Goyt~\cite{goy:apt} has considered the analogous question for
$\Pi_n(\cP)$ where $\cP\sbe\Pi_3$.  For example, 
\beq
\label{Fn}
\Pi_n(1/2/3,13/2)=F_n
\eeq
where $F_n$ is the $n$th Fibonacci number.

\subsection{Statistics}

The {\it inversion number\/} of  $p=a_1 a_2 \ldots a_n\in \fS_n$ is
$$
\inv p =\#\{(a_i,a_j)\ :\ \mbox{$i<j$ and $a_i>a_j$}\}.
$$
Also, the {\it major index\/} of $p$ is defined to be
$$
\maj p=\sum_{a_i>a_{i+1}} i.
$$
It  is well-known, and easy to prove, that if $q$ is a variable then
$$
\sum_{\pi\in\fS_n} q^{\inv\pi}=
\sum_{\pi\in\fS_n} q^{\maj\pi}=
1(1+q)(1+q+q^2)\cdots(1+q+\cdots+q^{n-1}).
$$
A statistic on $\fS_n$ with this generating function is said to be
{\it Mahonian\/} in honor of Major Percy MacMahon~\cite{mac:ca} 
who made the first systematic study of $\inv$ and $\maj$.  And the
polynomial product above is called a {\it $q$-analogue\/} of the
integer $n!$

Babson and Steingr\'{\i}msson~\cite{bs:gpp} defined generalized
permutation patterns by insisting that certain elements of the pattern
be adjacent in the larger permutation.  They then showed that most
Mahonian statistics in the literature can be written as linear
combinations of the statistics defined by generalized patterns.  

The Stirling numbers of the the second kind, $S(n,k)$, count 
the number of set partitions of
$[n]$ with $k$ blocks.
Carlitz~\cite{car:af,car:qbn} introduced a $q$-analogue, $S_q(n,k)$, of
$S(n,k)$.  Milne~\cite{mil:rgf}, Garsia and Remmel~\cite{gr:qcr},
Leroux~\cite{ler:rmq}, and Wachs and White~\cite{ww:pqs} have all
given set partition analogues of the inv statistic whose generating
function is $S_q(n,k)$, possibly up to a factor of $q^{{k\choose2}}$.
Sagan~\cite{sag:mss} and later White~\cite{whi:isp} gave maj
statistics for set partitions.   In the previously mentioned paper of
Goyt~\cite{goy:apt}, generalized patterns for set partitions are
defined.  He then uses them to obtain various statistics in the
literature as well as enumerates the number of partitions which avoid
them.

Carlitz~\cite{car:qfn} was also the first to define a $q$-analogue,
$F_n(q)$, of the Fibonacci numbers.  These polynomials and their
generalizations have
been extensively studied by Cigler~\cite{cig:qfpFQ,cig:ncq,cig:qfpAC}
as well as Shattuck and Wagner~\cite{sw:mfc}.  In view of~\ree{Fn},
one can define related $q$-analogues using the generating functions
for various set partition statistics over the family $\Pi_n(1/2/3,13/2)$.
This yields a new and unified approach to the study of $F_n(q)$ and
its relatives which is being pursued by Goyt and Sagan~\cite{gs:sps}.

\subsection{Partial orders}

The set $\fS$ of all permutations becomes a poset (partially
ordered set) by defining $p\le q$ if and only if there is a copy of
$p$ in $q$.  One of the fundamental invariants of any poset is its
M\"obius function, $\mu$.  See Stanley's 
text~\cite[Chapter 3]{sta:ec1} for information about posets in general
and the M\"obius function in particular.  Wilf asked the following
question.
\begin{question}[Wilf~\cite{wil:pp}]
If $p\le q$ in $\fS$ then what is $\mu(p,q)$?
\end{question}

This question has been partially answered as follows.  Call a
permutation $p$ {\it layered\/} if it has the form
$$
p=i,i-1,\ldots,1,i+j,i+j-1,\ldots,i+1,i+j+k,i+j+k-1,\ldots,i+j+1,\ldots
$$
where $i,j,k,\ldots$ are called the {\it layer lengths of $p$\/}.
There is a bijection between the layered permutations in $\fS_n$ and
compositions (ordered integer partitions) of $n$ gotten by sending $p$
as above to the composition $(i,j,k,\ldots)$.  Denoting the set of all
compositions by $\bbP^*$, we have a partial order on this set induced
by the pattern containment order on $\fS$.  This partial order was
first studied by Bergergon, Bousquet-M\'elou, and Dulucq~\cite{bbd:spc}
who counted its saturated lower chains.  Further work in this direction
was done by Snellman~\cite{sne:scc,sne:spa}.
The M\"obius function of $\bbP^*$ was obtained by Sagan and
Vatter~\cite{sv:mfc} in two ways, combinatorially and using discrete
Morse theory.  It was also rederived by Bj\"orner and
Sagan~\cite{bs:rmf} using the theory of regular languages.  This poset turns
out to be intimately related to subword order, whose M\"obius function
was first completely determined by Bj\"orner~\cite{bjo:mfs,bjo:mff}
and again by Bj\"orner and Reutenauer~\cite{br:rmf}.

Of course, we can partially order $\Pi$ by pattern containment and ask
the same question.
\begin{question}
If $\pi\le \si$ in $\Pi$ then what is $\mu(\pi,\si)$?
\end{question}
There is clearly a bijection between layered permutations and
layered partitions.  So the work cited above
applies to this poset as well.  Note that if we restrict the full poset of
compositions to the compositions which only contain ones and twos,
then we get a corresponding partial order on $\Pi(1/2/3,13/2)$ having
rank numbers equal to the $F_n$.
Goyt~[personal communication] is
currently investigating what can be said in various other posets
related to $\Pi$ whose rank numbers are given by certain generalized Fibonacci numbers. 

\bigskip

\bibliographystyle{acm}
\begin{small}
\bibliography{ref}
\end{small}

\end{document}